\documentclass[12pt]{amsart}

\usepackage[utf8]{inputenc} 
\usepackage[T1]{fontenc}    
\usepackage{lmodern}        
\usepackage{amssymb}        
\usepackage{amsthm}         
\usepackage{mathtools}      
\usepackage{mathrsfs}
\usepackage{physics}        
\usepackage{asymptote}      
\usepackage{subcaption}     
\usepackage{xcolor}         
\usepackage[
    hyperfootnotes=false      
]{hyperref}                 
\usepackage{cleveref}       
\usepackage{microtype}      
\usepackage{booktabs}       
\usepackage{tabularx}       
\usepackage[
    maxbibnames=99,           
    maxcitenames=2,           
    backend=biber             
]{biblatex}                 
\usepackage{filecontents}   
\usepackage{mathptmx}
\usepackage[margin=1.5in]{geometry}
\usepackage{float}
\usepackage{layouts}

\DeclareFieldFormat[article]{pages}{#1}
\DeclareFieldFormat[article]{volume}{\mkbibbold{#1}}
\makeatletter
\def\l@section{\@tocline{1}{0pt}{0pc}{5pc}{}}
\def\l@subsection{\@tocline{2}{0pt}{2.5pc}{5pc}{}}
\makeatother

\definecolor{linkblue}{HTML}{00356B}
\definecolor{linkgold}{HTML}{DA9100}
\definecolor{linkred}{RGB}{159, 29, 53}
\hypersetup{
        breaklinks = true, 
	colorlinks = true,
	linkcolor = linkred,
	citecolor = linkgold,
	urlcolor = linkblue    
}

\theoremstyle{plain}
\newtheorem{theorem}{Theorem}[section]
\crefname{theorem}{Theorem}{Theorems}

\crefname{conjecture}{Conjecture}{Conjectures}
\newtheorem{proposition}[theorem]{Proposition} 
\crefname{proposition}{Proposition}{Propositions}
 
\crefname{corollary}{Corollary}{Corollaries}
 
\crefname{lemma}{Lemma}{Lemmas}
\crefname{ineq}{inequality}{inequalities}
\newtheorem{remark}[theorem]{Remark}
\crefname{remark}{Remark}{Remarks}

\theoremstyle{definition}

\crefname{example}{Example}{Examples}
\newtheorem{definition}[theorem]{Definition}

\crefname{appendix}{Appendix}{Appendices}
\crefname{section}{Section}{Sections}
\crefname{figure}{Figure}{Figures}
\crefname{table}{Table}{Tables}

\newcommand{\field}[1]{\mathbf{#1}}
\newcommand{\R}{\field{R}}

\title[]{Numerically Computed Double and Triple Bubbles in $\R^3$ for Density $r^p$}

\author[]{Eve Parrott}

\address{Livingston High School, 30 Robert Harp Drive, Livingston, NJ 07039}
\email{evep1260@gmail.com}

\date{\today}

\begin{document}


\begin{asydef}
import geometry;
import patterns;

// Drawing Colors
pen backgroundcolor = white;

picture dot(real Hx=1mm, real Hy=1mm, pen p=currentpen)
{
  picture tiling;
  path g=(0,0)--(Hx,Hy);
  draw(tiling,g,invisible);
  dot(tiling, (0,0), p+linewidth(1));
  return tiling;
}

// Add fill patterns
add("hatch", hatch(1mm, currentpen));
add("dot",dot(1mm, currentpen));

pair bubble(real radius1, real radius2, pen drawingpen = currentpen, int lab = 0) {

	// Draws the appropriate double bubble for bubbles of volume
	// pi * radius1^2 and pi * radius2^2, respectively
	// Calculate the length of the segment connecting the centers of the
	// circles, from the law of cosines. Note that the angle opposite this
	// length formed by the radii is 60 degrees.
	real length = sqrt(radius1^2 + radius2^2 - radius1 * radius2);

	pair origin = (0,0);
	path circle1 = circle(origin, radius1);
	path circle2 = circle((length,0), radius2);

	// Check if the circles have equal curvature. If so, create a straight-line
	// path between their intersection points. Otherwise, draw the bulge as
	// an arc between the intersection points with curvature equal to the 
	// difference of the two circles.
	pair[] i_points = intersectionpoints(circle1, circle2);
	path bulge;
	
	if (radius1 == radius2) {
		// If the radii are equal, the "bulge" will actually be a circle 
		// through infinity, which mathematically is just a line but
		// computationally is an error. In this case, just draw the segment
		// connecting the two points.
		bulge = i_points[0] -- i_points[1];
	} else {
		
		// Calculate the radius of the middle circular arc
		real bulge_radius = 1 / abs(1/radius1 - 1/radius2);

		// Calculate the center of the middle circular arc
		// 	(1) Draw two circles of bulge_radius centered
		//	    at the two singularities. They will intersect
		//	    one another at the two possible centers.
		//	(2) Check which center is correct by comparing
		//	    the original radii.
		path big_circ_1 = circle(i_points[0], bulge_radius);
		path big_circ_2 = circle(i_points[1], bulge_radius);
		pair[] big_i_points = intersectionpoints(big_circ_1, big_circ_2);
		pair big_center;
		if (radius1 > radius2) {
			big_center = big_i_points[1];
		} else {
			big_center = big_i_points[0];
		}
		
		// Calculate the angles from the horizontal between the center of
		// the bulge and the two singularities. Draw an arc at the center
		// using those angles. Since arcsin(...) has a limited range,
		// we must add 180 degrees to theta if the center is on the
		// right, along with reversing the path direction.
		real height = arclength(i_points[0] -- i_points[1]) / 2.0;
		real theta = aSin(height/bulge_radius);
		
		if (radius1 > radius2) {
			bulge = arc(big_center, bulge_radius, -theta + 180, theta + 180);
		} else {
			bulge = arc(big_center, bulge_radius, theta, -theta);
		}
	}

	draw(bulge, drawingpen);

	// Calculate the parameterized time when the og circles intersect
	real[][] i_times = intersections(circle1, circle2);
	real time_1_a = i_times[0][0]; // Time on circle1 of first intersection 
	real time_1_b = i_times[1][0]; // Time on circle1 of second intersection 
	real time_2_a = i_times[0][1]; // Time on circle2 of first intersection 
	real time_2_b = i_times[1][1]; // Time on circle2 of second intersection 

	// The second path is weird to account for the fact that paths always 
	// start at time 0, and we can't run the path in reverse. This then 
	// calculates how long it takes to walk the whole path, and then breaks 
	// the arc up into subpaths.
	// path bubble1 = subpath(circle1, time_1_a, time_1_b) -- bulge -- cycle;
	path bubble1 = subpath(circle1, time_1_a, time_1_b) -- reverse(bulge) -- cycle;
	path bubble2 = subpath(circle2, 0, time_2_a) -- bulge -- subpath(circle2, time_2_b, arctime(circle2, arclength(circle2))) -- cycle;

	filldraw(bubble1, pattern("hatch"), drawingpen);
	filldraw(bubble2, pattern("dot"), drawingpen);

	if (lab == 1) {
		label(Label("$\bm{\Omega_1}$",Fill(backgroundcolor)),origin);
		label(Label("$\bm{\Omega_2}$",Fill(backgroundcolor)),(length, 0));
	}
	
	return i_points[0];
}
\end{asydef}

\begin{abstract}
Using Brakke's Evolver, we numerically verify previous conjectures for optimal double bubbles for density $r^p$ in $R^3$ and our own new conjectures for triple bubbles.
\end{abstract}

\maketitle


\section{Introduction}

\begin{figure}[h!]
    \begin{subfigure}[t]{0.49\linewidth}
        \centering
        \includegraphics[width=0.8\linewidth]{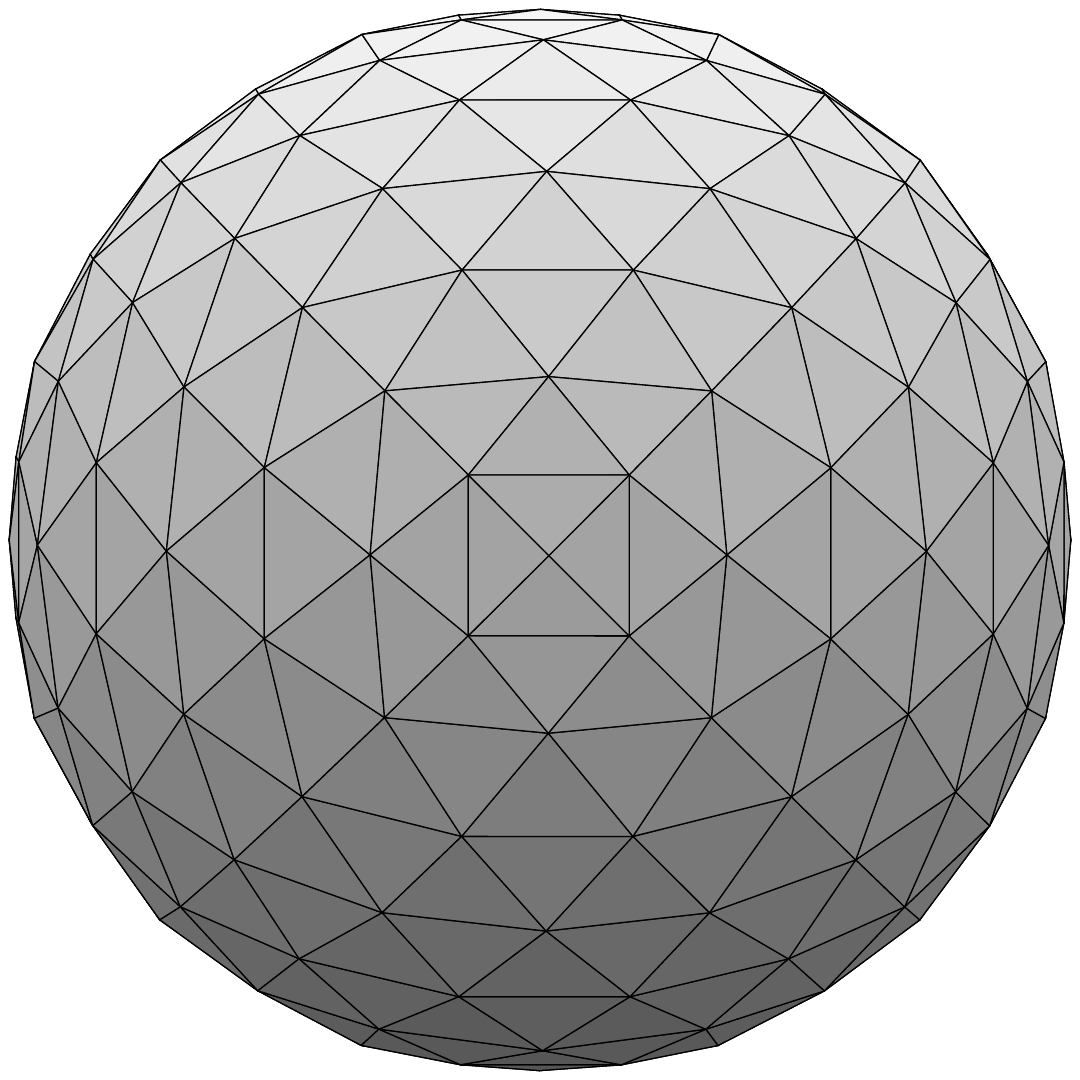}
    \end{subfigure}
    \hfill 
    \begin{subfigure}[t]{0.49\linewidth}
        \centering
        \includegraphics[width=\linewidth]{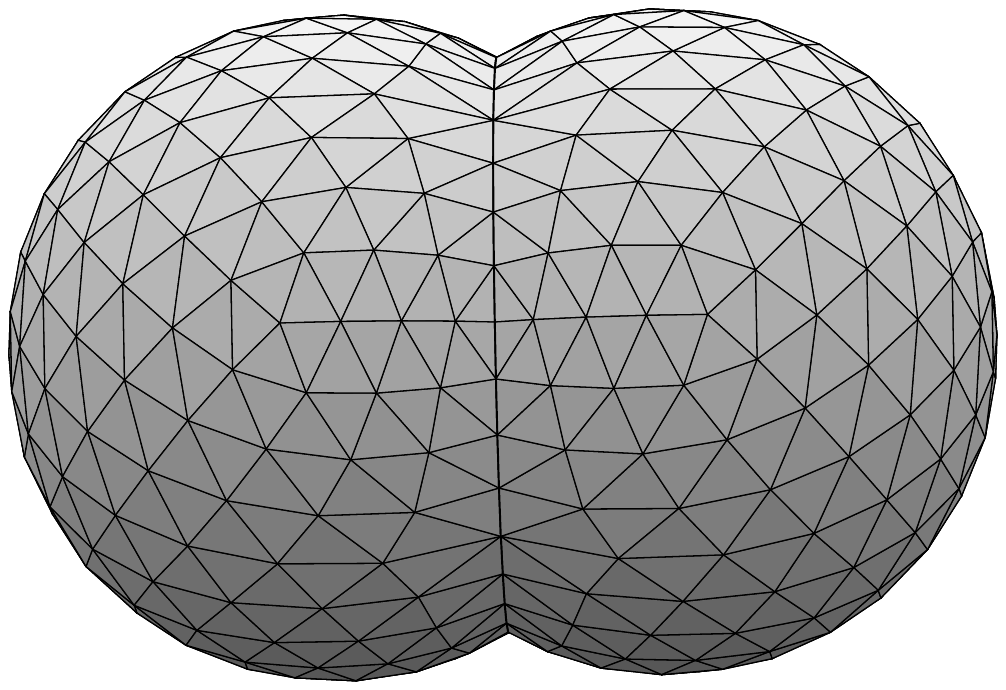}
        \begin{minipage}{.1cm}
        \vfill 
            \end{minipage} 
    \end{subfigure}
    \newline
    \begin{subfigure}[t]{0.49\linewidth}
        \centering
        \includegraphics[width=0.8\linewidth]{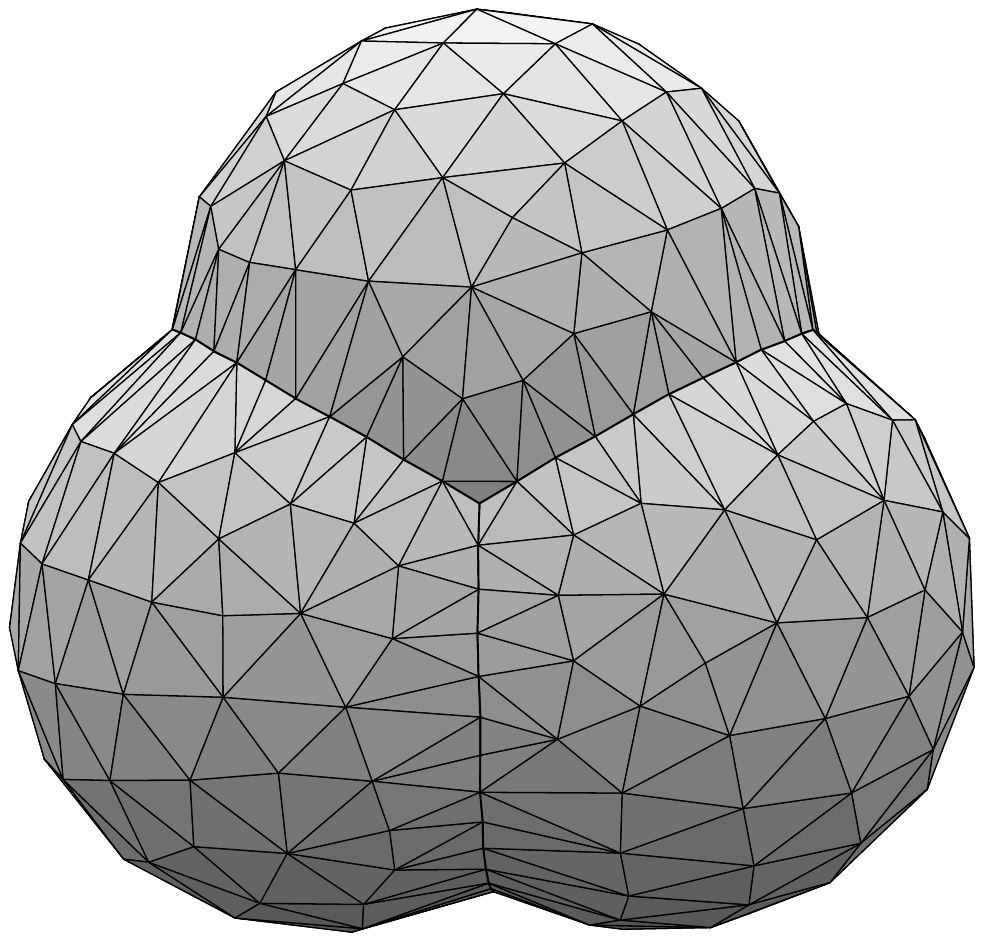}
    \end{subfigure}
    \hfill
    \caption{The optimal single, double \cite{annals}, and triple bubbles \cite{Milman2022} with equal volumes in Euclidean space.}
    \label{fig:euclidean-bubbles}
\end{figure}

The isoperimetric problem is one of the oldest in mathematics. 
It asks for the least-perimeter way to enclose given volume. For a single volume in Euclidean space (with uniform density) of any dimension, the well-known solution is any sphere. With 
density $r^p$, \textcite{G14} found that the solution for a single volume is a sphere \emph{through} the origin. 
For \emph{two} volumes in
Euclidean space, \textcite{annals}
showed that the standard double bubble, consisting of three spherical caps meeting along a sphere in threes at $120^\circ$ angles as in \cref{fig:euclidean-bubbles}, provides the 
isoperimetric cluster. Recently, \textcite{Milman2022} showed that the standard triple bubble is the isoperimetric cluster (see \cref{fig:euclidean-bubbles}). 

In this paper we use Brakke's Evolver \cite{Brakke2013} to numerically compute double and triple bubbles in $\R^3$ with density $r^p$ for various $p > 0$, extending many of the results of Collins \cite{Collins2023} from $\R^2$ to $\R^3$. Some videos are downloadable from \href{https://drive.google.com/drive/folders/1-3AJLlG08R-lC1MB2gCfqxSAUm48EPOx?usp=sharing}{Google Drive}.  \cref{doubleprop} supports the conjecture of Hirsch et al. \cite{Morgan2021} that the optimal double bubble is the standard double bubble with a singular circle passing through the origin (\cref{fig:double-bubbles}). \cref{tripleprop} indicates that the optimal triple bubble resembles a standard triple bubble (\cref{fig:euclidean-bubbles}) with one vertex at the origin (\cref{fig:triple-bubbles}).
\subsection*{History} Examination of isoperimetric regions in the plane with 
density $r^p$ began in \citeyear{G06} when \textcite{G06} showed that
the isoperimetric solution for a single area in the plane with density $r^p$ is 
a convex set containing the origin. It was something of a surprise when 
\textcite{dahlberg2010} proved that the solution is a circle
through the origin.
In \citeyear{G14} \textcite{G14} extended this result to higher dimensions. In \citeyear{china19} \textcite{china19} studied the $1$-dimensional case,
showing that the best single bubble is an interval with one endpoint at the 
origin and that the best double bubble is a pair of adjacent intervals which meet at the origin. Ross \cite{Ross2022} showed that in $\R^1$ multiple bubbles start with the two smallest meeting at the origin and the rest in increasing order alternating side to side.

\subsection*{Acknowledgements}

This paper was written under the guidance of Frank Morgan remotely, following separate work done at the MathPath summer camp. The author would like to thank Professor Morgan for his advice and support. Additionally, the author would like to thank Professor Ken Brakke for generous help with his Evolver. Finally, the author would like to thank Marcus Collins for contributing his research and insight.

\section{Definitions}

\begin{definition}[Density Function]
Given a smooth Riemannian manifold $M$, a \emph{density} on $M$ is a positive 
continuous function (perhaps vanishing at isolated points) that weights each point $p$ in $M$ with a certain mass 
$f(p)$. Given a region $\Omega \subset M$ with piecewise smooth boundary, the weighted volume (or area) and boundary measure or perimeter of $\Omega$ are given by
\[
V(\Omega) = \int_{\Omega} f \dd{V_0} \quad\text{and}\quad P(\Omega) = 
\int_{\partial \Omega} f \dd{P_0},
 \]
where $\dd{V_0}$ and $\dd{P_0}$ denote Euclidean volume and perimeter. We may also refer to the perimeter of $\Omega$ as the perimeter of its boundary.
\end{definition}

\begin{definition}[Isoperimetric Region]
A region $\Omega \subset M$ is 
\emph{isoperimetric} if it has the smallest weighted perimeter of all regions 
with the same weighted volume. The boundary of an isoperimetric region is also 
called isoperimetric.
\end{definition}

We can generalize the idea of an isoperimetric region by considering two or more
volumes.
\begin{definition}[Isoperimetric Cluster]
An isoperimetric cluster is a set of disjoint open regions~$\Omega_i$ of 
volume~$V_i$ such that 
the perimeter of the union of the boundaries is minimized.
\end{definition}

For example, in $R^3$ with density 1, optimal clusters are known for one volume, two volumes (\textcite{annals}), three volumes (\textcite{Milman2022}), as in \cref{fig:euclidean-bubbles}. Note that for density $r^p$, scalings of minimizers are minimizers, because scaling up by a factor of $\lambda$ scales perimeter by $\lambda^{p+2}$ and volume by $\lambda^{p+3}$. \\

The following proposition summarizes existence, boundedness, and regularity of isoperimetric clusters in $\R^n$ with density, proved by Hirsch et al. \cite{Morgan2021} following such results for single bubbles (Rosales et al. \cite[Thm.~2.5]{rosales}) and Morgan \cite{morgan}.

\begin{proposition}[Existence, Boundedness, and Regularity]\label{existence-boundedness-regularity}
    Consider $\R^n$ with radial non-decreasing density f, smooth and positive except possibly at the origin, such that $f(r) \to \infty$ as $r \to \infty.$ An isoperimetric cluster that encloses and separates given volumes exists (Hirsch et al.  \cite[Thm.~2.5]{Morgan2021} ) and is bounded (Hirsch et al. \cite[Prop.~2.6]{Morgan2021}). In $\R^3$, the cluster consists of smooth constant-generalized-mean-curvature surfaces meeting in threes at $120^\circ$ angles along curves, which in turn meet in fours at $\cos^{-1}(-\frac{1}{3}) \approx 109^\circ$ (\textcite[Section 13.9]{morgan}), except possibly at the origin.
\end{proposition}

\begin{definition}[Generalized Curvature]\label{def:gen_curve}
We define, as in \cite{rosales} and \cite[Chap. 3]{DeRosa2023} the (generalized) mean curvature of a surface with respect to the inward-pointing unit normal $N$ as the function
\[
\mathrm{H_{\psi}} = nH - \langle \nabla \psi, \mathbf{N} \rangle,
\]
where $H$ is the classical mean curvature of a surface (the sum of the principal curvatures) and $\psi$ is the logarithm of the density function $f$. 
This comes from the first variation of perimeter formula, so that generalized curvature has the interpretation as minus the rate at which $\log f$ changes in the direction perpendicular to the curve (see \textcite[Sect.\ 3]{rosales}).
\end{definition}

\begin{remark}
    In $\R^3$ with density $r^p$, a circular arc has constant generalized curvature if and only if the sphere passes through or is centered at the origin  \cite[Rem.~2.9]{Morgan2021}.
    \label{genconstantcurvatureremark}
\end{remark} 
\section{Double and Triple Bubbles in \texorpdfstring{$\R^3$}{r^2}  with density \texorpdfstring{$r^p$}{rp}}


With Brakke's Evolver \cite{Brakke2013}, \cref{doubleprop} supports the conjecture of Hirsch et al. \cite{Morgan2021} that the optimal double bubble in $\R^3$ with density $r^p$ is a standard double bubble with the singular circle passing through the origin. \cref{tripleprop} provides a conjecture on the form of the triple bubble. 

\begin{figure}[H]
    \centering
    \begin{subfigure}[c]{0.3\textwidth}
        \centering
        \includegraphics[width=\linewidth]{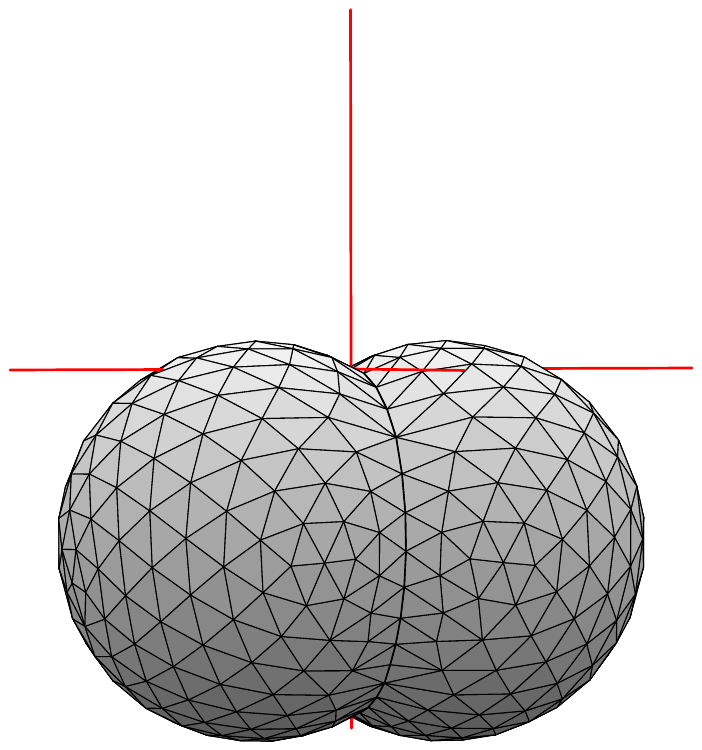 }
        \caption{Equal volumes}
        \label{fig:equal-double-bubble}
    \end{subfigure}
    \hfill 
    \begin{subfigure}[c]{0.3\textwidth}
        \centering
        \includegraphics[width=\linewidth]{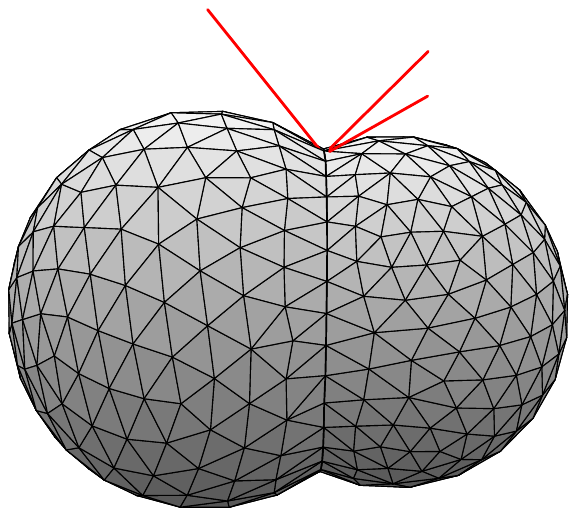}
        \caption{Volumes $8$ and $4$}
        \label{fig:unequal-double-bubble}
    \end{subfigure}
    \hfill 
    \begin{subfigure}[c]{0.3\textwidth}
        \centering
        \includegraphics[width=1.2\linewidth]{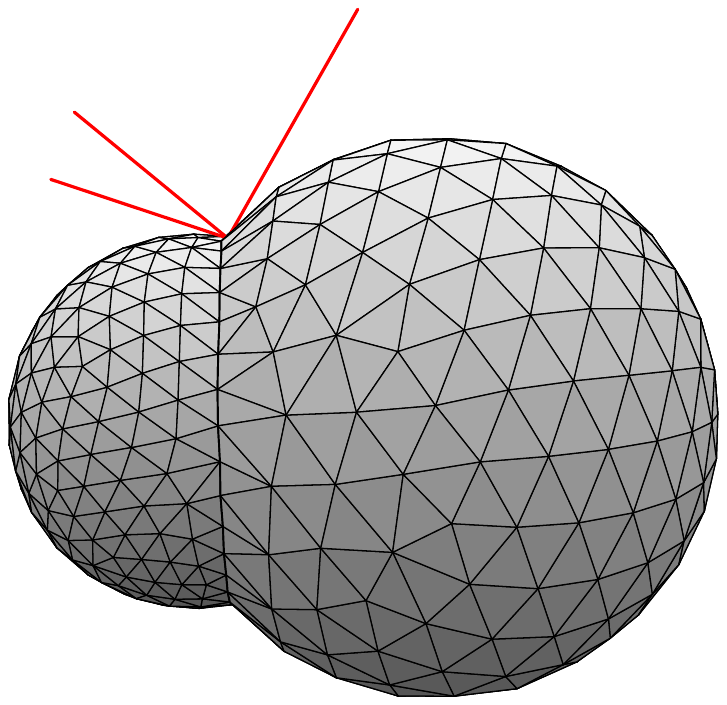}
        \caption{Volumes $0.8$ and $8$}
        \label{fig:double-bubble-3}
        \begin{minipage}{.1cm}
        \vfill
        \end{minipage}
    \end{subfigure}
    \newline
    \begin{subfigure}[c]{0.3\textwidth}
        \centering
        \includegraphics[width=\linewidth]{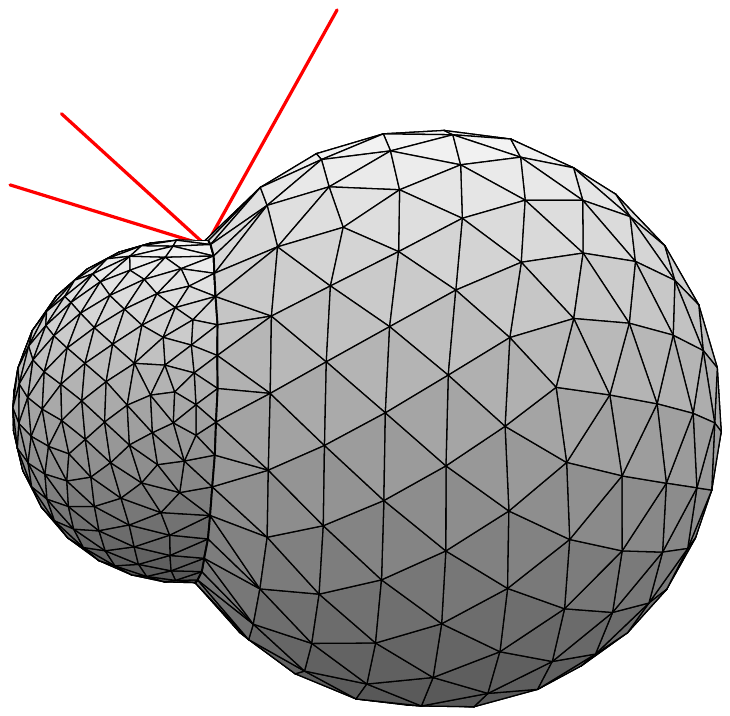}
        \caption{Volumes $0.4$ and $8$}
        \label{fig:double-bubble-4}
        \begin{minipage}{.1cm}
        \vfill
        \end{minipage}
    \end{subfigure}
    \hfill
    \begin{subfigure}[c]{0.3\textwidth}
        \centering
        \includegraphics[width=\linewidth]{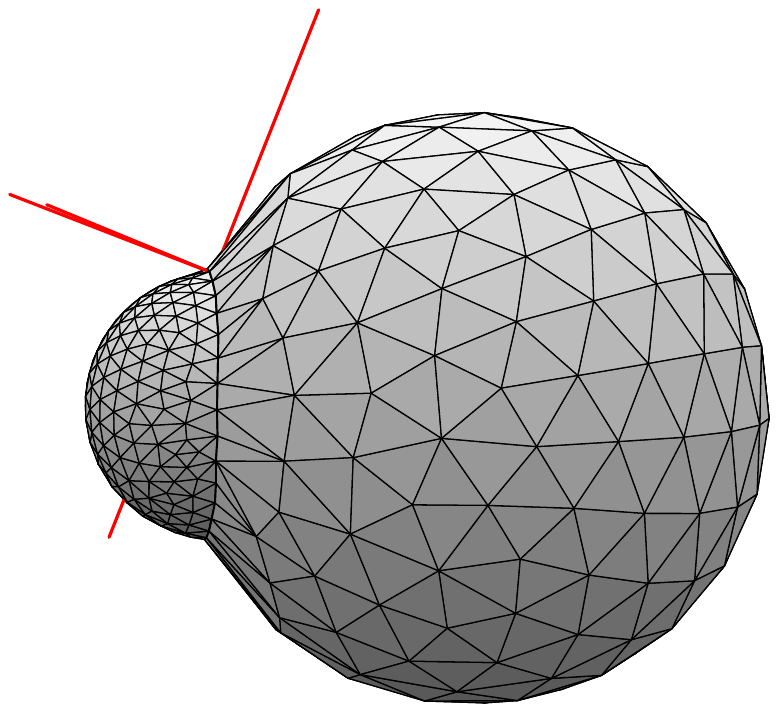}
        \caption{Volumes $0.08$ and $8$}
        \label{fig:double-bubble-5}
    \end{subfigure}
    \hfill
    \begin{subfigure}[c]{0.3\textwidth}
        \centering
        \includegraphics[width=\linewidth]{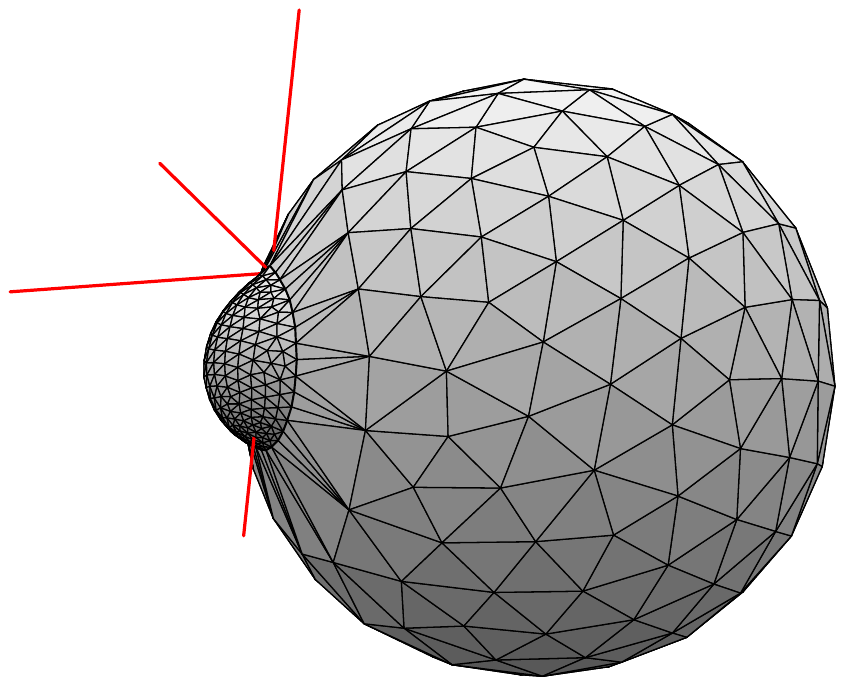}
        \caption{Volumes $0.008$ and $8$}
        \label{fig:double-bubble-6}
        \begin{minipage}{.1cm}
        \vfill
        \end{minipage}
    \end{subfigure}
    \caption{Computations in Brakke's Evolver \cite{Brakke2013} in $\R^3$ with density $r^2$ support the conjecture that the optimal double bubble is a standard double bubble with the singular circle passing through the origin (where the red axes meet). Densities $r^5, r^3$, and $r^{0.5}$ are apparently identical.}
    \label{fig:double-bubbles}
\end{figure}

\begin{proposition}[Double Bubble] \label{doubleprop} 
Computations on Brakke's Evolver \cite{Brakke2013} support the conjecture \cite{Morgan2021} that the optimal double bubble in $\R^3$ for density $r^p$ consists of a standard double bubble with the singular circle passing through the origin. See \cref{fig:double-bubbles}.
\end{proposition}

\begin{figure}[H]
    \centering
    \begin{subfigure}[c]{0.3\textwidth}
        \centering
        \includegraphics[width=\linewidth]{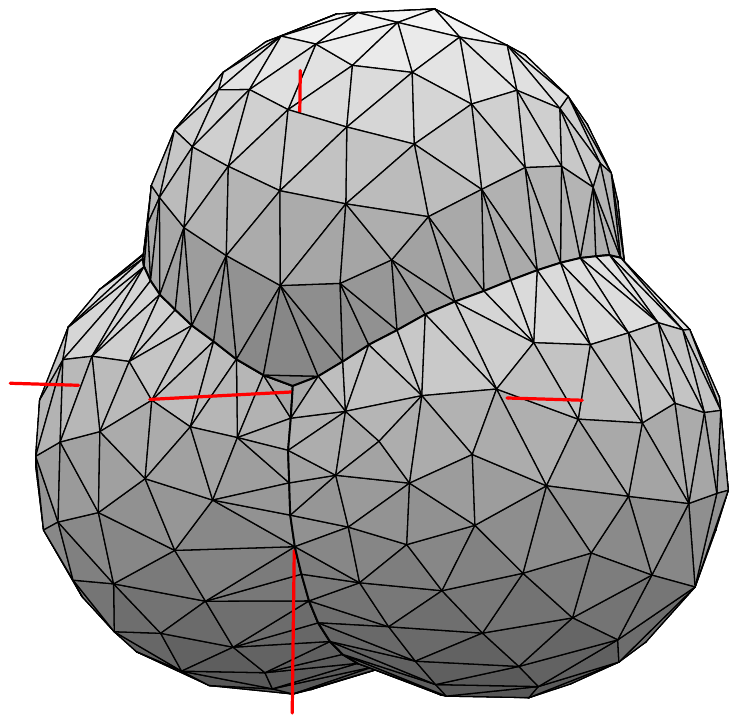}
        \caption{Equal volumes of $300$}
        \label{fig:triple-bubble-1}
    \end{subfigure}
    \hfill
    \begin{subfigure}[c]{0.3\textwidth}
        \centering
        \includegraphics[width=\linewidth]{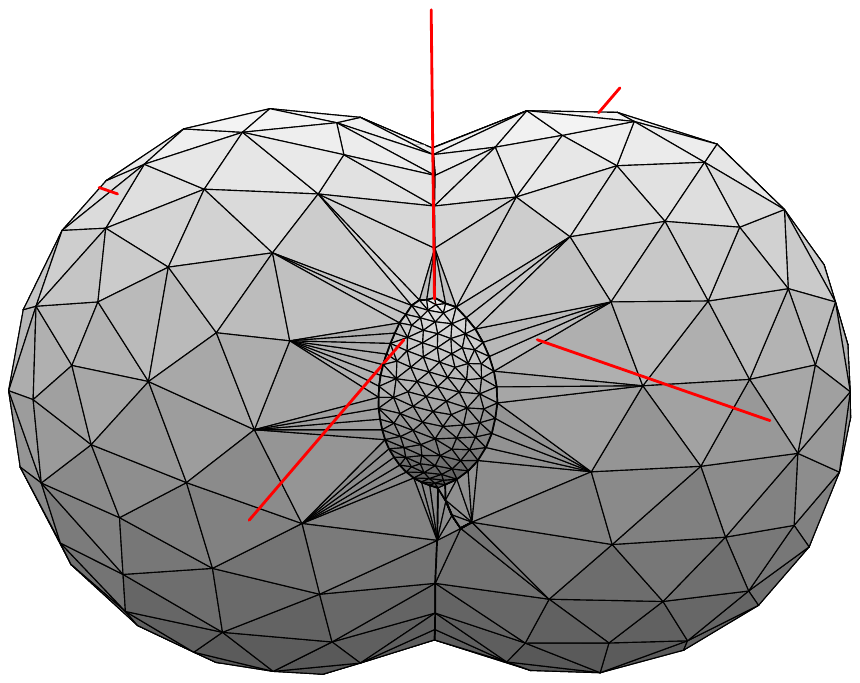}
        \caption{Volumes $0.3$, $300$, and $300$}
        \label{fig:triple-bubble-5}
    \end{subfigure}
    \hfill 
    \begin{subfigure}[c]{0.3\textwidth}
        \centering
        \includegraphics[width=\linewidth]{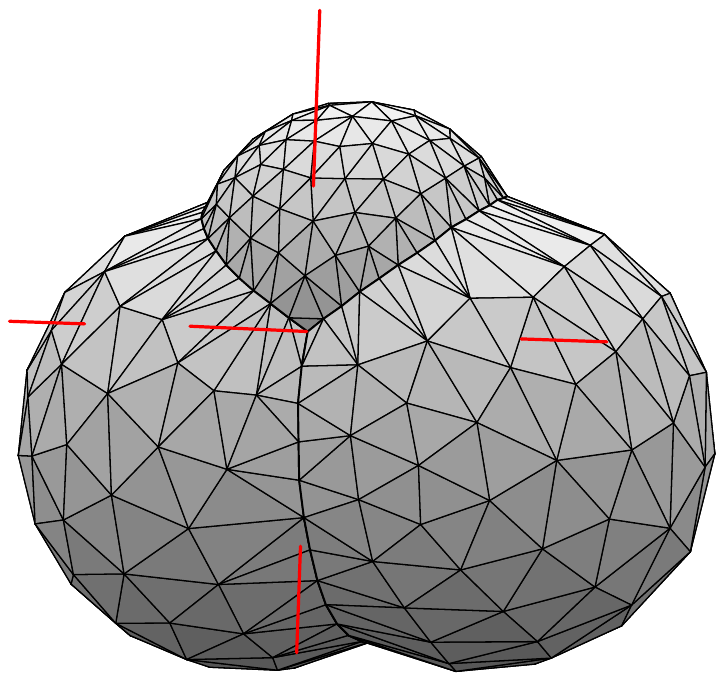}
        \caption{Volumes $300$, $300$, and $30$}
        \label{fig:triple-bubble-3}
        \begin{minipage}{1.cm}
        \vfill
        \end{minipage}
    \end{subfigure}
    \newline 
    \begin{subfigure}[c]{0.3\textwidth}
        \centering
        \includegraphics[width=\linewidth]{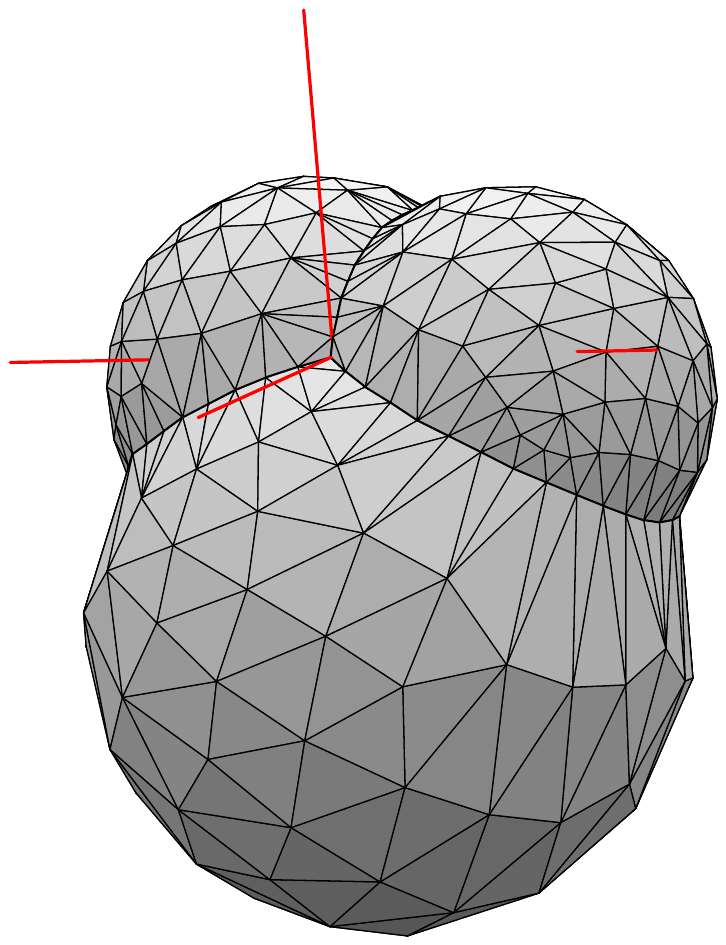}
        \caption{Volumes $0.3$, $0.3$, and $3$}
        \label{fig:triple-bubble-grey}
    \end{subfigure}
    \hfill 
    \begin{subfigure}[c]{0.3\textwidth}
        \centering
        \includegraphics[width=\linewidth]{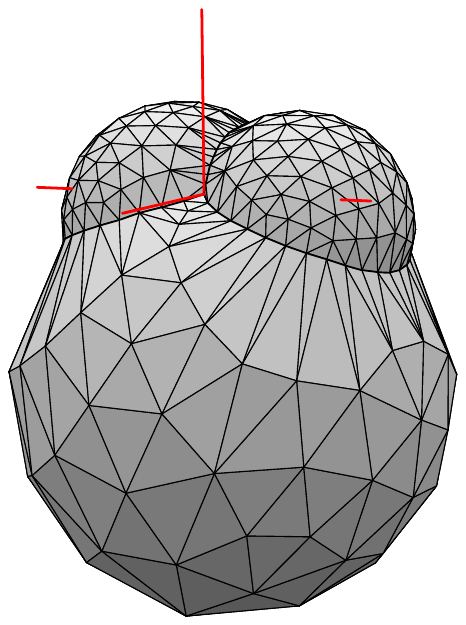}
        \caption{Volumes $0.3$, $0.3$, and $15$}
        \label{fig:triple-bubble-green}
    \end{subfigure} 
    \newline
    \begin{subfigure}[c]{0.3\textwidth}
        \centering
        \includegraphics[width=\linewidth]{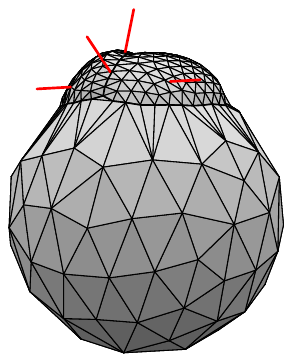}
        \caption{Volume $0.3$, $0.3$, and $60$}
        \label{fig:triple-bubble-blue}
    \end{subfigure}
    \hfill
    \begin{subfigure}[c]{0.3\textwidth}
        \centering
        \includegraphics[width=\linewidth]{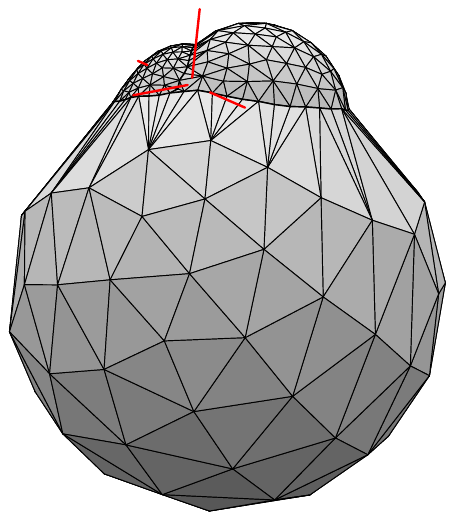}
        \caption{Volumes $1$, $5$, and $1000$}
        \label{fig:triple-bubble-6}
    \end{subfigure}
    \hfill
    \begin{subfigure}[c]{0.3\textwidth}
        \centering
        \includegraphics[width=1.2\linewidth]{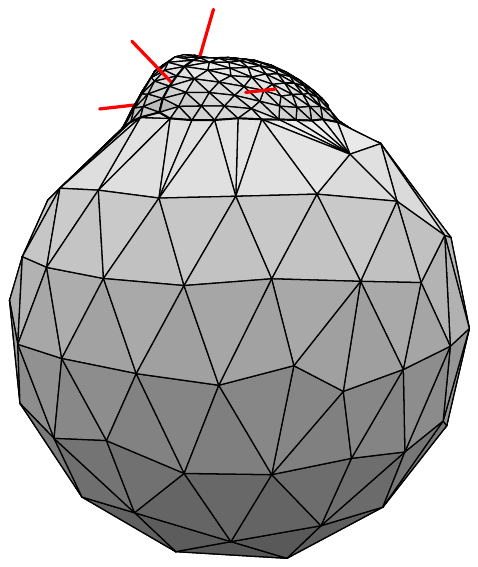}
        \caption{Volumes $0.3$, $0.3$, and $300$}
        \label{fig:triple-bubble-red}
    \end{subfigure}
    \hfill
    \caption{Computations in Brakke's Evolver \cite{Brakke2013} in $\R^3$ with density $r^2$ suggest that the optimal triple bubble resembles a standard triple bubble with one vertex at the origin.}
    \label{fig:triple-bubbles}
\end{figure}

\begin{proposition}[Triple Bubble]
\label{tripleprop}
Computations with Brakke's Evolver \cite{Brakke2013} indicate that the optimal triple bubble in $\R^3$ with density $r^p$ resembles a standard triple bubble with one vertex at the origin. See \cref{fig:triple-bubbles}.
\end{proposition}

\begin{figure}[h!]
    \centering
    \begin{subfigure}[c]{0.49\textwidth}
        \centering
        \includegraphics[width=.4\linewidth]{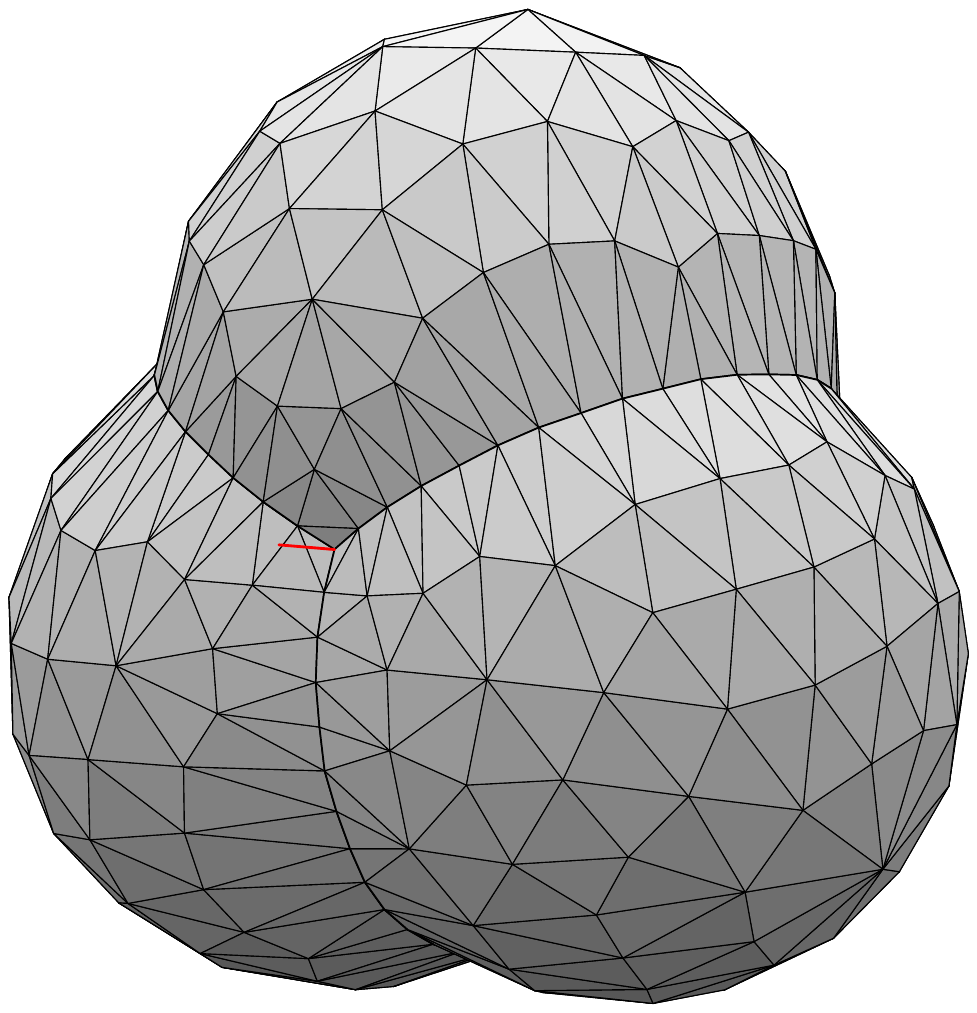}
        \label{fig:0.5-triple}
    \end{subfigure}
    \hfill 
    \begin{subfigure}[c]{0.49\textwidth}
        \centering
        \includegraphics[width=.7\linewidth]{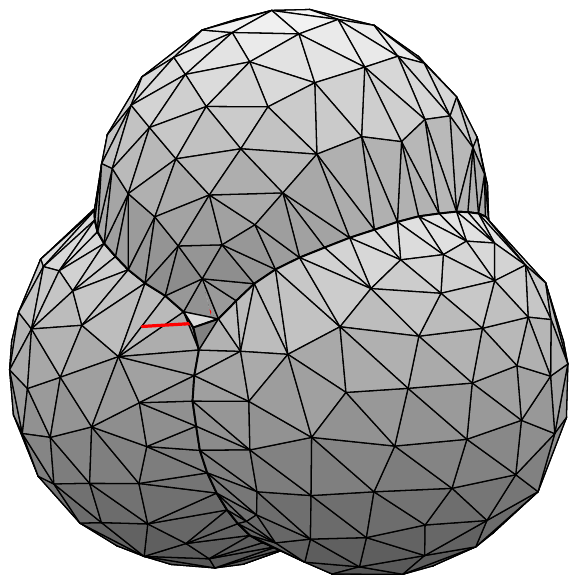}
        \label{fig:3-triple}
        \begin{minipage}{0.1cm}
        \vfill
        \end{minipage} 
    \end{subfigure}
    \caption{Computations in Brakke’s Evolver \cite{Brakke2013} in $R^3$ with densities $r^{0.5}$ and $r^3$ suggest that the optimal triple bubble resembles a standard triple bubble with one vertex at the origin, as for density $r^2$ of Figure 3.}
    \label{fig:triple-bubble-p<2}
\end{figure} 
\begin{remark}
    The phenomenon of vertices of triple bubbles in $R^2$ with density $r^p$ approaching as $p$ increases (Collins \cite[Prop. 3.2]{Collins2023}) apparently does not hold in $R^3$.
\end{remark}

\begin{figure}[h!]
    \centering
    \begin{subfigure}[c]{0.49\textwidth}
        \centering
        \includegraphics[width=0.65\linewidth]{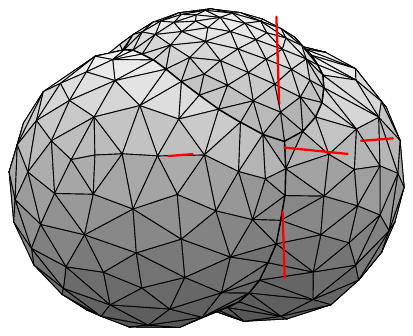}
        \caption{Our conjectured triple bubble with volumes $10$, $10$, and $1$ has surface area just over $32$.}
        \label{fig:triple-bubble-perimeter}
    \end{subfigure}
    \hfill 
    \begin{subfigure}[c]{0.49\textwidth}
        \centering
        \includegraphics[width=0.65\linewidth]{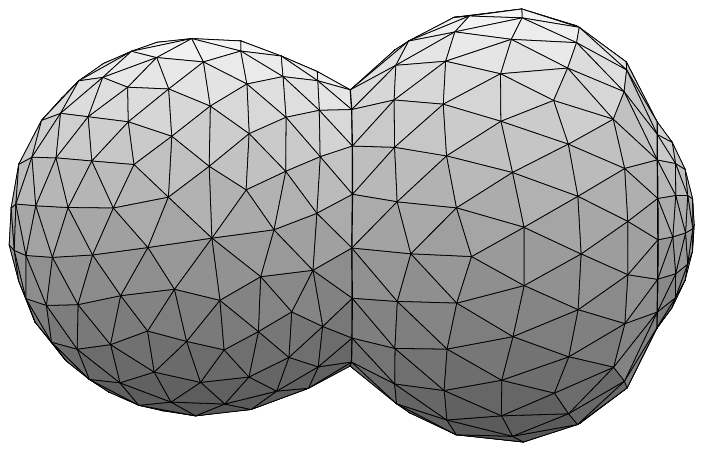}
        \caption{A linear chain with volumes $10$, $10$, and $1$ has surface area just over $36$.}
        \label{fig:triple-chain-perimeter}
        \begin{minipage}{0.1cm}
        \vfill
        \end{minipage}
    \end{subfigure}
    \caption{Our conjectured triple bubble has less surface area than a linear chain in $R^3$ with density $r^2$. Densities $r^{0.5}$ and $r^3$ are apparently similar.}
    \label{fig:triple-comparison}
\end{figure}

\begin{proposition}
\label{triplechainprop}
Our conjectured triple bubble of \cref{tripleprop} has a lower surface area than three bubbles in a linear chain, as in \cref{fig:triple-comparison}. 
\end{proposition}

\clearpage 
\newpage
\printbibliography
\end{document}